\documentclass[11pt]{amsart}

\usepackage{graphicx,amsfonts}
\usepackage{amsmath}
\usepackage{color}
\usepackage{txfonts}
\usepackage[colorlinks]{hyperref}
\usepackage{mathptmx}
\usepackage{mathpazo}
\usepackage{pifont}
\usepackage{graphicx}
\usepackage{amssymb}
\usepackage{tikz}
\usetikzlibrary{positioning}

\newcommand{\n}{\noindent}
\newtheorem{thm}{Theorem}

\newtheorem{cor}{Corollary}
\newtheorem{ex}{Example}
\begin{document}

\markboth{Nilanjan De, Sk. Md. Abu Nayeem, Anita Pal}
{F-index of some graph operations}

%
%

\title{F-Index of Some Graph Operations}

\author{Nilanjan De}

\address{Department of
Basic Sciences and Humanities (Mathematics),
Calcutta Institute of Engineering and Management, Kolkata, India.}
\email{de.nilanjan@rediffmail.com}

\author{Sk. Md. Abu Nayeem}

\address{Department of Mathematics, Aliah University,
IIA/27, New Town, Kolkata - 700 156, India.}
\email{nayeem.math@aliah.ac.in}

\author{Anita Pal}

\address{Department of
Mathematics,
National Institute of Technology, Durgapur, India.}
\email{anita.buie@gmail.com}

\maketitle

\begin{abstract}
The F-index of a graph is defined as the sum of cubes of the vertex degrees of the graph. This was introduced in 1972, in the same paper where the first and second Zagreb indices were introduced to study the structure-dependency of total  $\pi$-electron energy. But this topological index was not further studied till then. Very recently, Furtula and Gutman [B. Furtula, I. Gutman, A forgotten topological index, \textit{J. Math. Chem.}, \textbf{53(4)}(2015) 1184--1190.] reinvestigated the index and named it ``forgotten topological index" or ``F-index". In that paper, they present some basic properties of this index and showed that this index can enhance the physico-chemical applicability of Zagreb index. Here, we study the behavior of this index under several graph operations and apply our results to find the F-index of different chemically interesting molecular graphs and nano-structures.\\[5pt]
\end{abstract}

\section{Introduction}

Suppose $G$ is a simple connected graph and $V(G)$ and $E(G)$ denote the vertex set and edge set of $G$, respectively.  For any vertex ${v}\in V(G)$, let ${{d}_{G}}(v)$ denote its degree, that is the number of neighbors of $v$ and $N(v)$ denote the set of vertices which are the neighbors of the vertex $v$, so that $|N(v)|={{d}_{G}}(v)$. In chemistry, biochemistry and nanotechnology different topological indices are found to be useful in isomer discrimination, structure-property relationship, structure-activity relationship and pharmaceutical drug design.
The first and second Zagreb indices of a graph $G$, denoted by $M_1(G)$ and $M_2(G)$, are among the oldest, most popular and extremely studied vertex-degree based topological indices and are respectively defined as
\[{{M}_{1}}(G)=\sum\limits_{v\in V(G)}{{{d}_{G}}{{(v)}^{2}}}=\sum\limits_{uv\in E(G)}{[{{d}_{G}}(u)+{{d}_{G}}(v)]}\] and \[{{M}_{2}}(G)=\sum\limits_{uv\in E(G)}{{{d}_{G}}(u){{d}_{G}}(v)}.\]

These indices were introduced in a paper in 1972 \cite{gutm72} to study the structure-dependency of the total $\pi$-electron energy ($\varepsilon$). It was found that the $\varepsilon$ depends on $M_1(G)$ and thus provides a measure of carbon skeleton of the underlying molecules. In the same paper, another topological index, defined as sum of cubes of degrees of the vertices of the graph was also shown to influence $\varepsilon$. However this index was not further studied till then, except in a recent article by Furtula and Gutman \cite{fur15} where they reinvestigated this index and studied some basic properties of this index. They showed that the predictive ability of this index is almost similar to that of first Zagreb index and for the entropy and acentric factor, both of them yield correlation coefficients greater than 0.95. They named this index as  ``forgotten topological index" or ``F-index". Throughout the present paper we name this index as F-index and denote it by $F(G)$, so that
\[F(G)=\sum\limits_{v\in V(G)}{{{d}_{G}}{{(v)}^{3}}}=\sum\limits_{uv\in E(G)}{[{{d}_{G}}{{(u)}^{2}}+{{d}_{G}}{{(v)}^{2}}]}.\]

As we know that some chemically interesting graphs can be obtained by different graph operations on some general or particular graphs, it is important to study such graph operations in order to understand how it is related to the corresponding topological indices of the original graphs. In \cite{kha09}, Khalifeh et al. derived some exact formulae for computing first and second Zagreb indices under some graph operations. In \cite{das13}, Das et al. derived some upper bounds for multiplicative Zagreb indices for different graph operations. In \cite{aza14}, Azari presented some lower bounds for Narumi-Katayama index under several graph operations. In \cite{de14}, the present authors computed some bounds and exact formulae of the connective eccentric index under different graph operations. There are several other results regarding various topological indices under different graph operations are available in the literature. In \cite{aza13}, Azari and Iranmanesh presented explicit formulas for computing the eccentric-distance sum of different graph operations. Interested readers are referred to \cite{Ash10,kha08,tava14,veyl15,esk13,aza13a,aza15,nd15} in this regard.

In this paper, we present some exact expressions for the F-index of different graph operations such as union, join, Cartesian product, composition, tensor product, strong product, corona product, generalized hierarchical product, disjunction, symmetric difference, splice and link of graphs. Also we apply our results to compute the F-index for some important classes of molecular graphs and nano-structures.

\section{Main Results and Discussions}
In this section, we study F-index of various graph operations like union, join, Cartesian product, composition, tensor product, strong product, corona product, generalized hierarchical product, disjunction, symmetric difference, link and splice of graphs. These operations are binary and all the graphs are connected, finite and simple. In the following, if not indicated otherwise, we use the notation $V(G_i)$ for the vertex set, ${{E}(G_i)}$ for the edge set, ${{n}_{i}}$ for the number of vertices and ${{m}_{i}}$ for the number of edges of the graph ${{G}_{i}}$, $i\in \left\{ 1,2,\ldots,k \right\}$, respectively. Throughout the paper, we use the familiar notations $P_n$, $C_n$ and $K_n$ to denote a path graph, cycle graph and complete graph with $n$ number of vertices, respectively.

\subsection{Union}
Let $G_1$, $G_2$,...,$G_k$ be $k$ graphs with disjoint vertex sets. Then their union ${{G}_{1}}\cup {{G}_{2}}\cup ...\cup {{G}_{k}}$ is the graph with vertex set ${V({G}_{1})}\cup {V({G}_{2})}\cup ...\cup {V({G}_{k})}$ and the edge set ${E({G}_{1})}\cup {E({G}_{2})}\cup ...\cup {E({G}_{k})}$. The degree of a vertex $v$ of ${{G}_{1}}\cup {{G}_{2}}\cup ...\cup {{G}_{k}}$ is equal to the degree of the vertex $v$ in the component $G_i$, $i=1,2,...,k$, that contains it.
In the following theorem we obtain the F-index of the union of $k$ number of graphs.
\begin{thm}The F-index of ${{G}_{1}}\cup {{G}_{2}}\cup ...\cup {{G}_{k}}$ is given by
\[F({{G}_{1}}\cup {{G}_{2}}\cup ...\cup {{G}_{k}})=F({{G}_{1}})+F({{G}_{2}})+...+F({{G}_{k}}).\]
\end{thm}
\n\textit{Proof.} By definition of F-index we have
\begin{eqnarray*}
F({{G}_{1}}\cup {{G}_{2}}\cup ...\cup {{G}_{k}})&=&\sum\limits_{v\in V({{G}_{1}}\cup {{G}_{2}}\cup ...\cup {{G}_{k}})}{({{d}_{{{G}_{1}}\cup {{G}_{2}}\cup ...\cup {{G}_{k}}}}(v)})^3\\
&=&\sum\limits_{v\in V({{G}_{1}})}{{{d}_{{{G}_{1}}}}{{(v)}^{3}}}+\sum\limits_{v\in V({{G}_{2}})}{{{d}_{{{G}_{2}}}}{{(v)}^{3}}}+...+\sum\limits_{v\in V({{G}_{k}})}{{{d}_{{{G}_{k}}}}{{(v)}^{3}}}\\
&=&F({{G}_{1}})+F({{G}_{2}})+...+F({{G}_{k}}),
\end{eqnarray*}
which completes the proof. \qed

\subsection{Join}
The join ${{G}_{1}}+{{G}_{2}}$ of two graphs $G_1$ and $G_2$ is the union ${{G}_{1}}\cup{{G}_{2}}$ together with all the edges joining ${V({G}_{1})}$ and ${V({G}_{2})}$. The degree of a vertex $v$ of ${{G}_{1}}+{{G}_{2}}$ is
\[{{d}_{{{G}_{1}}+{{G}_{2}}}}(v) = \left\{ \begin{array}{ll}
{{d}_{{{G}_{1}}}}(v)+{{n}_{2}},v\in V({{G}_{1}})\\[2mm]
{{d}_{{{G}_{2}}}}(v)+{{n}_{1}},v\in V({{G}_{2}}).
\end{array}\right.\]
In general, for $k$ graphs $G_1$, $G_2$,...,$G_k$, the degree of a vertex $v$ in ${{G}_{1}}+{{G}_{2}}+...+{{G}_{k}}$ is given by ${d}_{{{G}_{1}}+{{G}_{2}}+...+{G}_{k}}(v)$ $= {{d}_{{{G}_{i}}}}(v)+{n}-{{n}_{i}}$, where $v$ is originally a vertex of the graph $G_i$ and $n={{n}_{1}}+{{n}_{2}}+...+{{n}_{k}}$.

In the following theorem we compute the F-index of the join of $k$ number of graphs.
\begin{thm} The F-index of ${{G}_{1}}+{{G}_{2}}+...+{{G}_{k}}$ is given by
\[F({{G}_{1}}+{{G}_{2}}+...+{{G}_{k}})=\sum\limits_{i=1}^{k}{F({{G}_{i}})}+3\sum\limits_{i=1}^{k}{{{{\bar{n}}}_{i}}{{M}_{1}}({{G}_{i}})}+6\sum\limits_{i=1}^{k}{{{{\bar{n}}}_{i}}^{2}{{m}_{i}}}+\sum\limits_{i=1}^{k}{{{n}_{i}}{{{\bar{n}}}_{i}}^{3}},\]
where ${{\bar{n}}_{i}}=n-{{n}_{i}}, i=1,2,...,k$ and $n={{n}_{1}}+{{n}_{2}}+...+{{n}_{k}}$.
\end{thm}

\n\textit{Proof.} We have
\begin{eqnarray*}
F({{G}_{1}}+{{G}_{2}}+...+{{G}_{k}})&=&\sum\limits_{i=1}^{k}{\sum\limits_{v\in V({{G}_{i}})}{{{({{d}_{{{G}_{i}}}}(v)+{{{\bar{n}}}_{i}})}^{3}}}}\\
                          &=&\sum\limits_{i=1}^{k}{\sum\limits_{v\in V({{G}_{i}})}{({{d}_{{{G}_{i}}}}{{(v)}^{3}}+3{{{\bar{n}}}_{i}}{{d}_{{{G}_{i}}}}{{(v)}^{2}}+3{{{\bar{n}}}_{i}}^{2}{{d}_{{{G}_{i}}}}(v)}+{{{\bar{n}}}_{i}}^{3}})\\
                             &=&\sum\limits_{i=1}^{k}{\sum\limits_{v\in V({{G}_{i}})}{{{d}_{{{G}_{i}}}}{{(v)}^{3}}}}+3\sum\limits_{i=1}^{k}{{{{\bar{n}}}_{i}}\sum\limits_{v\in V({{G}_{i}})}{{{d}_{{{G}_{i}}}}{{(v)}^{2}}}}\\
                             &&+3\sum\limits_{i=1}^{k}{{{{\bar{n}}}_{i}}^{2}\sum\limits_{v\in V({{G}_{i}})}{{{d}_{{{G}_{i}}}}(v)}}+\sum\limits_{i=1}^{k}{{{n}_{i}}{{{\bar{n}}}_{i}}^{3}},
\end{eqnarray*}
which completes the proof. \qed

Let ${{G}_{1}}={{G}_{2}}=...={{G}_{p}}=G$ and $pG$ denote the join of $p$ copies of $G$. Then the following corollaries follow as direct consequence of the previous theorem.
\begin{cor}
Let, $n$ and $m$ be the number of vertices and edges of $G$, respectively. Then
\[F(pG)=pF(G)+3np(p-1){{M}_{1}}(G)+6{{n}^{2}}mp{{(p-1)}^{2}}+{{n}^{4}}p{{(p-1)}^{3}}.\]
\end{cor}

\begin{cor}
The F-index of ${{G}_{1}}+{{G}_{2}}$ is given by
\[F({{G}_{1}}+{{G}_{2}})=F({{G}_{1}})+F({{G}_{2}})+3{{n}_{2}}{{M}_{1}}({{G}_{1}})+3{{n}_{1}}{{M}_{1}}({{G}_{2}})+6{{n}_{2}}^{2}{{m}_{1}}+6{{n}_{1}}^{2}{{m}_{2}}+{{n}_{1}}{{n}_{2}}^{3}+{{n}_{2}}{{n}_{1}}^{3}.\]
\end{cor}
The suspension of a graph $G$ is defined as ${{K}_{1}}+G$. So from the Corollary 2 the following result follows.
\begin{cor}
Let, $n$ and $m$ be the number of vertices and edges of $G$, respectively. Then the F-index of suspension of $G$ is given by
\[F({{K}_{1}}+G)=F(G)+3{{M}_{1}}(G)+{{n}^{3}}+6m+n.\]
\end{cor}

\begin{ex}
The complete $n$-partite graph ${{K}_{{{m}_{1}},{{m}_{2}},...,{{m}_{n}}}}$ (Fig.1) on ${{m}_{1}}+{{m}_{2}}+...+{{m}_{n}}$ vertices can be considered as ${{\bar{K}}_{{{m}_{1}}}}+{{\bar{K}}_{{{m}_{2}}}}+...+{{\bar{K}}_{{{m}_{n}}}}$. Then the F-index of ${{K}_{{{m}_{1}},{{m}_{2}},...,{{m}_{n}}}}$ is given by \[F({{K}_{{{m}_{1}},{{m}_{2}},...,{{m}_{n}}}})=\sum\limits_{i=1}^{n}{{{m}_{i}}{{{\bar{m}}}_{i}}^{3}},\] where ${{\bar{m}}_{i}}=({{m}_{1}}+{{m}_{2}}+...+{{m}_{n}})-{{m}_{i}}, i=1,2,...,n$.
\end{ex}

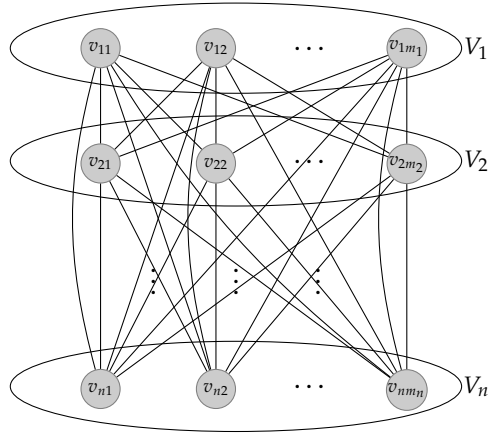
\begin{figure}[h]
\begin{center}
\begin{tikzpicture}[place/.style={circle,draw=black!50,fill=black!20,
inner sep=0pt, minimum size = 5.2mm}, bend angle=15]
\draw (1.8,0) ellipse (30mm and 6mm);
\draw (1.8,-1.5) ellipse (30mm and 6mm);
\draw (1.8,-4.5) ellipse (30mm and 6mm);
\draw (2.8,0) node{\ldots};
\draw (2.8,-1.5) node{\ldots};
\draw (1.8,-3) node{\vdots \hspace{1cm}\vdots \hspace{1cm}\vdots};
\draw (2.8,-4.5) node{\ldots};
\draw (5,0) node{\footnotesize $V_1$};
\draw (5,-1.5) node{\footnotesize $V_2$};
\draw (5,-4.5) node{\footnotesize $V_n$};
{\tiny
\node[place] (1) {$v_{11}$};
\node[place] (2) [right=of 1] {$v_{12}$};
\node[place] (3) [right=of 2, xshift=1cm] {$v_{1m_1}$};
\node[place] (4) [below=of 1] {$v_{21}$}
edge (1) edge (2) edge (3);
\node[place] (5) [below=of 2] {$v_{22}$}
edge (1) edge (2) edge (3);
\node[place] (6) [below=of 3] {$v_{2m_2}$}
edge (1) edge (2) edge (3);
\node[place] (7) [below=of 1, yshift=-3cm]{$v_{n1}$}
edge[bend left] (1) edge (2) edge (3) edge (4) edge (5) edge (6);
\node[place] (8) [below=of 2, yshift=-3cm] {$v_{n2}$}
edge (1) edge[bend left] (2) edge (3) edge (4) edge (5) edge (6);
\node[place] (9) [below=of 3, yshift=-3cm] {$v_{nm_n}$}
edge[bend angle=10, bend left] (1) edge (2) edge[bend left] (3) edge (4) edge (5) edge (6);
}
\end{tikzpicture}
\caption{\label{f3}The complete $n$-partite graph.}
\end{center}
\end{figure}

\begin{ex}
The wheel graph ${{W}_{n}}$ on $(n+1)$ vertices is the suspension of ${{C}_{n}}$ and the fan graph ${{F}_{n}}$ on $(n+1)$ vertices is the suspension of ${{P}_{n}}$. So their F-indices are given by

(i) $F({{W}_{n}})={{n}^{3}}+27n$,

(ii) $F({{F}_{n}})={{n}^{3}}+27n-38$.
\end{ex}

\begin{ex}
The dutch windmill graph or flower graph is the suspension of $m$ copies of $K_2$, denoted by $m{K_2}$. So its F-index is given by $F({K_1}+m{K_2})=8m^3+16m$.
\end{ex}
\begin{ex}
The cone graph ${C}_{m,n}$ is defined as ${C_m}+\bar{K}_n$. So its F-index is calculated as $F({C}_{m,n})=m{n^3}+{m^3}n+6m{n^2}+12mn+8m$.
\end{ex}

\subsection{Cartesian Product}

The Cartesian product of $G_1$ and $G_2$, denoted by $G_1\otimes G_2$, is the graph with vertex set $V(G_1)\times V(G_2)$ and any two vertices $({{u}_{p}},{{v}_{r}})$ and $({{u}_{q}},{{v}_{s}})$ are adjacent if and only if [${{u}_{p}}={{u}_{q}}$ and ${{v}_{r}}{{v}_{s}}\in E(G_2)$] or [${{v}_{r}}={{v}_{s}}$ and ${{u}_{p}}{{u}_{q}}\in E(G_1)$]. In the following theorem first we find the first Zagreb index of the Cartesian product of $k$ number of graphs.

\begin{thm}Let $n$ be the total number of vertices in $\bigotimes_{i=1}^{k} G_i$, then the first Zagreb index of ${{G}_{1}}\otimes {{G}_{2}}\otimes ...\otimes {{G}_{k}}$ is given by
\[{M_1}\left(\bigotimes_{i=1}^k G_i\right)= {n}\sum\limits_{i=1}^{k}{\frac{{{M}_{1}}({{G}_{i}})}{{{n}_{i}}}+}4{n}\sum\limits_{\substack{i,j=1\\i\ne j}}^{k}{\frac{{{m}_{i}}}{{{n}_{i}}}\cdot\frac{{{m}_{j}}}{{{n}_{j}}}}.\]
\end{thm}
\n\textit{Proof.} For the proof of this theorem we refer to Theorem 1 of \cite{kha09}.\qed

In the following theorem we obtain the F-index of the Cartesian product of $k$ number of graphs.
\begin{thm}
The F-index of ${{G}_{1}}\otimes {{G}_{2}}\otimes ...\otimes {{G}_{k}}$ is given by
\[F({{G}_{1}}\otimes {{G}_{2}}\otimes ...\otimes {{G}_{k}})=n\sum\limits_{i=1}^{k}{\frac{F({{G}_{i}})}{{{n}_{i}}}}+6n\sum\limits_{\substack{i,j=1\\i\ne j}}^{k}{\frac{{{M}_{1}}({{G}_{i}})}{{{n}_{i}}}\cdot\frac{{{m}_{j}}}{{{n}_{j}}}}+8n\sum\limits_{\substack{p,q,r=1\\p\ne q\ne r}}^{k}{\frac{{{m}_{p}}}{{{n}_{p}}}\cdot\frac{{{m}_{q}}}{{{n}_{q}}}}\cdot\frac{{{m}_{r}}}{{{n}_{r}}}.\]
\end{thm}

\n\textit{Proof.} First we prove the result for $k=2$.

We have, ${{d}_{{{G}_{1}}\otimes {{G}_{2}}}}(a,b)={{d}_{{{G}_{1}}}}(a)+{{d}_{{{G}_{2}}}}(b)$.
So, from definition of F-index we have
\begin{eqnarray*}
F({{G}_{1}}\otimes {{G}_{2}})&=&{{\sum\limits_{a\in V({{G}_{1}})}{\sum\limits_{b\in V({{G}_{2}})}{({{d}_{{{G}_{1}}}}(a)+{{d}_{{{G}_{2}}}}(b))}}}^{3}}\\
&=&\sum\limits_{b\in V({{G}_{2}})}{\sum\limits_{a\in V({{G}_{1}})}{{{d}_{{{G}_{1}}}}{{(a)}^{3}}}}+\sum\limits_{a\in V({{G}_{1}})}{\sum\limits_{b\in V({{G}_{2}})}{{{d}_{{{G}_{2}}}}{{(b)}^{3}}}}\\
&&+3\sum\limits_{a\in V({{G}_{1}})}{{{d}_{{{G}_{1}}}}{{(a)}^{2}}\sum\limits_{b\in V({{G}_{2}})}{{{d}_{{{G}_{2}}}}(b)+}}3\sum\limits_{a\in V({{G}_{1}})}{{{d}_{{{G}_{1}}}}(a)\sum\limits_{b\in V({{G}_{2}})}{{{d}_{{{G}_{2}}}}{{(b)}^{2}}}}\\
&=&{{n}_{2}}F({{G}_{1}})+{{n}_{1}}F({{G}_{2}})+6{{m}_{2}}{{M}_{1}}({{G}_{1}})+6{{m}_{1}}{{M}_{1}}({{G}_{2}}).
\end{eqnarray*}

Let $n'$ be the total number of vertices in $\bigotimes_{i=1}^{k-1} G_i$. Then by an inductive argument, using Theorem 3, we have
\begin{eqnarray*}
\displaystyle F\left(\bigotimes_{i=1}^k G_i\right)&=&F\left(\bigotimes_{i=1}^{k-1} G_i\bigotimes G_k\right)\\
&=&n_kF\left(\bigotimes_{i=1}^{k-1} G_i\right)+n'F(G_k)+6m_kM_1\left(\bigotimes_{i=1}^{k-1} G_i\right)+6n'\sum_{i=1}^{k-1}\frac{m_i}{n_i}M_1(G_k)\\
&=&{{n}_{k}}\left[ {n}'\sum\limits_{i=1}^{k-1}{\frac{F({{G}_{i}})}{{{n}_{i}}}+}6{n}'\sum\limits_{\substack{i,j=1\\i\ne j}}^{k-1}{\frac{{{M}_{1}}({{G}_{i}})}{{{n}_{i}}}\cdot\frac{{{m}_{j}}}{{{n}_{j}}}+8{n}'\sum\limits_{\substack{p,q,r=1\\p\ne q\ne r}}^{k-1}{\frac{{{m}_{p}}}{{{n}_{p}}}\cdot\frac{{{m}_{q}}}{{{n}_{q}}}\cdot\frac{{{m}_{r}}}{{{n}_{r}}}}} \right]\\
&&+{n}'F({{G}_{k}})+6{{m}_{k}}\left[ {n}'\sum\limits_{i=1}^{k-1}{\frac{{{M}_{1}}({{G}_{i}})}{{{n}_{i}}}+}4{n}'\sum\limits_{\substack{i,j=1\\i\ne j}}^{k-1}{\frac{{{m}_{i}}}{{{n}_{i}}}\cdot\frac{{{m}_{j}}}{{{n}_{j}}}} \right]+6{n}'{{M}_{1}}({{G}_{k}})\sum\limits_{i=1}^{k-1}{\frac{{{m}_{i}}}{{{n}_{i}}}}\\
&=&\left[ n\sum\limits_{i=1}^{k-1}{\frac{F({{G}_{i}})}{{{n}_{i}}}+}n\frac{F({{G}_{k}})}{{{n}_{k}}} \right]+\left[ 6n\sum\limits_{\substack{i,j=1\\i\ne j}}^{k-1}\frac{{{M}_{1}}({{G}_{i}})}{{{n}_{i}}}\cdot\frac{{{m}_{j}}}{{{n}_{j}}}+6n\frac{{{m}_{k}}}{{{n}_{k}}}\sum\limits_{i=1}^{k-1}{\frac{{{M}_{1}}({{G}_{i}})}{{{n}_{i}}}}\right.\\
&&\left.+6n\frac{{{M}_{1}}({{G}_{k}})}{{{n}_{k}}}\sum\limits_{i=1}^{k-1}{\frac{{{m}_{i}}}{{{n}_{i}}}} \right]+\left[ 8n\sum\limits_{\substack{p,q,r=1\\p\ne q\ne r}}^{k-1}{\frac{{{m}_{p}}}{{{n}_{p}}}\cdot\frac{{{m}_{q}}}{{{n}_{q}}}\cdot\frac{{{m}_{r}}}{{{n}_{r}}}}+24n\frac{{{m}_{k}}}{{{n}_{k}}}\sum\limits_{\substack{i,j=1\\i\ne j}}^{k-1}{\frac{{{m}_{i}}}{{{n}_{i}}}\cdot\frac{{{m}_{j}}}{{{n}_{j}}}} \right],
\end{eqnarray*}

from where the desired result follows. \qed

\begin{cor}
The F-index of of the torus ${{C}_{{{n}_{1}}}}\otimes {{C}_{{{n}_{2}}}}\otimes...\otimes {{C}_{{{n}_{k}}}}$ is given by

$F({{C}_{{{n}_{1}}}}\otimes {{C}_{{{n}_{2}}}}\otimes ...\otimes {{C}_{{{n}_{k}}}})=8{k^3}{{n}_{1}}{{n}_{2}}\ldots{{n}_{k}}.$

\end{cor}

\begin{cor}
The F-index of the Hamming graph ${{K}_{{{n}_{1}}}}\otimes {{K}_{{{n}_{2}}}}\otimes ...\otimes {{K}_{{{n}_{k}}}}$ is given by
\begin{eqnarray*}
F({{K}_{{{n}_{1}}}}\otimes {{K}_{{{n}_{2}}}}\otimes ...\otimes {{K}_{{{n}_{k}}}})&=&{{n}_{1}}{{n}_{2}}\ldots{{n}_{k}}(\sum\limits_{i=1}^{k}{{{({{n}_{i}}-1)}^{3}}}+3\sum\limits_{\substack{i,j=1\\i\ne j}}^{k}{{{({{n}_{i}}-1)}^{2}}({{n}_{j}}-1)}\\
&&+\sum\limits_{\substack{p,q,r=1\\p\ne q\ne r}}^{k}{({{n}_{p}}-1)({{n}_{q}}-1)({{n}_{r}}-1)})\\
&=&{({n_1}+{n_2}+...+{n_k}-k)^3}{n_1}{n_2}...{n_k}.
\end{eqnarray*}

\end{cor}
\begin{ex}
For $k$-dimensional hypercube ${{Q}_{k}}={{K}_{2}}\otimes {{K}_{2}}\otimes ...\otimes {{K}_{2}}$ (k times), by our calculation, we have $F({{Q}_{k}})={{2}^{k}}{{k}^{3}}$.
\end{ex}
\begin{ex}
The F-index of ${{K}_{{{n}_{1}}}}\otimes {{K}_{{{n}_{2}}}}$ torus is given by \[F({{K}_{{{n}_{1}}}}\otimes {{K}_{{{n}_{2}}}})={{n}_{1}}{{n}_{2}}{({n_1}+{n_2}-2 )^3}\]
\end{ex}

\begin{ex}
Let $R$ and $S$ denote a $C_4$ nanotube and nanotorus, respectively. Then $R\cong{P_n}\otimes{C_m}$ and $S\cong{C_n}\otimes{C_m}$, for some integers $n$ and $m$. Then by our calculation, $F(R)=64{m}{n}-74{m}$ and  $F(S)=64{m}{n}$. Also if, $T\cong{P_n}\otimes{P_m}$, then $F(T)=64{m}{n}-74{m}-74{n}+72$.
\end{ex}
\subsection{Composition}

The composition or lexicographic product of two graphs ${{G}_{1}}$ and ${{G}_{2}}$ is denoted by ${{G}_{1}}[{{G}_{2}}]$. The vertex set of ${{G}_{1}}[{{G}_{2}}]$ is $V({{G}_{1}})\times V({{G}_{2}})$ and the degree of a vertex $(a,b)$ of ${{G}_{1}}[{{G}_{2}}]$ is given by ${{d}_{{{G}_{1}}[{{G}_{2}}]}}(a,b)={{n}_{2}}{{d}_{{{G}_{1}}}}(a)+{{d}_{{{G}_{2}}}}(b)$ and any two vertices $({{u}_{1}},{{u}_{2}})$ and $({{v}_{1}},{{v}_{2}})$ are adjacent if and only if ${{u}_{1}}{{v}_{1}}\in E({{G}_{1}})$ or [${{u}_{1}}={{v}_{1}}$ and  ${{u}_{2}}{{v}_{2}}\in E({{G}_{2}})$]. In the following theorem we compute the F-index of the composition of two graphs.
\begin{thm} The F-index of ${{G}_{1}}[{{G}_{2}}]$ is given by
\begin{eqnarray*}
F({{G}_{1}}[{{G}_{2}}])={{n}_{2}}^{4}F({{G}_{1}})+{{n}_{1}}F({{G}_{2}})+6{{n}_{2}}^{2}{{m}_{2}}{{M}_{1}}({{G}_{1}})+6{{n}_{2}}{{m}_{1}}{{M}_{1}}({{G}_{2}}).
\end{eqnarray*}
\end{thm}

\n\textit{Proof.} From definition of F-index, we have
\begin{eqnarray*}
F({{G}_{1}}[{{G}_{2}}])&=&\sum\limits_{(a,b)\in V({{G}_{1}}[{{G}_{2}}])}{{{d}_{{{G}_{1}}[{{G}_{2}}]}}}{{(a,b)}^{3}}\\
                 &=&\sum\limits_{a\in V({{G}_{1}})}{\sum\limits_{b\in V({{G}_{2}})}{{{[{{n}_{2}}{{d}_{{{G}_{1}}}}(a)+{{d}_{{{G}_{2}}}}(b)]}^{3}}}}\\
                 &=&{{n}_{2}}^{3}\sum\limits_{b\in V({{G}_{2}})}{\sum\limits_{a\in V({{G}_{1}})}{{{d}_{{{G}_{1}}}}{{(a)}^{3}}}}            +\sum\limits_{a\in V({{G}_{1}})}{\sum\limits_{b\in V({{G}_{2}})}{{{d}_{{{G}_{2}}}}{{(b)}^{3}}}}\\
                 &&+3{{n}_{2}}^{2}\sum\limits_{a\in V({{G}_{1}})}{{{d}_{{{G}_{1}}}}{{(a)}^{2}}\sum\limits_{b\in V({{G}_{2}})}{{{d}_{{{G}_{2}}}}(b)+}}3{{n}_{2}}\sum\limits_{a\in V({{G}_{1}})}{{{d}_{{{G}_{1}}}}(a)\sum\limits_{b\in V({{G}_{2}})}{{{d}_{{{G}_{2}}}}{{(b)}^{2}}}}\\
                 &=&{{n}_{2}}^{4}F({{G}_{1}})+{{n}_{1}}F({{G}_{2}})+6{{n}_{2}}^{2}{{m}_{2}}{{M}_{1}}({{G}_{1}})+6{{n}_{2}}{{m}_{1}}{{M}_{1}}({{G}_{2}}),
                 \end{eqnarray*}
which completes proof.                                             \qed

\begin{ex}
The fence graph is the composition of ${{P}_{n}}$ and ${{P}_{2}}$ and the closed fence graph is the composition of ${{C}_{n}}$ and ${{P}_{2}}$. So from the previous theorem, we have

(i) $F({{P}_{n}}[{{P}_{2}}])=250n-392$,

(ii) $F({{C}_{n}}[{{P}_{2}}])=250n.$
\end{ex}
\subsection{Tensor Product}

The tensor product or Kronecker product of two graphs ${{G}_{1}}$ and ${{G}_{2}}$ is denoted by ${{G}_{1}}\times {{G}_{2}}$ and any two vertices $({{u}_{1}},{{v}_{1}})$ and $({{u}_{2}},{{v}_{2}})$ are adjacent if and only if ${{u}_{1}}{{u}_{2}}\in E({{G}_{1}})$ and ${{v}_{1}}{{v}_{2}}\in E({{G}_{2}})$.  The degree of a vertex $(a,b)$ of ${{G}_{1}}\times {{G}_{2}}$ is given by
${{d}_{{{G}_{1}}\times {{G}_{2}}}}(a,b)={{d}_{{{G}_{1}}}}(a){{d}_{{{G}_{2}}}}(b)$.
In the following theorem, the F-index of the tensor product of two graphs is computed.
\begin{thm}

The F-index of ${{G}_{1}}\times {{G}_{2}}$ is given by
$F({{G}_{1}}\times {{G}_{2}})=F({{G}_{1}})F({{G}_{2}}).$
\end{thm}

\n\textit{Proof.} From definition of F-index, we have

\[F({{G}_{1}}\times {{G}_{2}})=\sum\limits_{(a,b)\in V({{G}_{1}}\times {{G}_{2}})}{{{d}_{{{G}_{1}}\times {{G}_{2}}}}}{{(a,b)}^{3}}=\sum\limits_{a\in V({{G}_{1}})}{\sum\limits_{b\in V({{G}_{2}})}{{{[{{d}_{{{G}_{1}}}}(a){{d}_{{{G}_{2}}}}(b)]}^{3}}}}=F({{G}_{1}})F({{G}_{2}}).\]
\qed

\begin{ex}
(i) $F({{P}_{n}}\times {{P}_{m}})=(8n-14)(8m-14)$,

(ii) $F({{C}_{n}}\times {{C}_{m}})=64nm$,

(iii) $F({{K}_{n}}\times {{K}_{m}})=nm{{(n-1)}^{3}}{{(m-1)}^{3}}$,

(iv) $F({{P}_{n}}\times {{C}_{m}})=8m(8n-14)$,

(v)  $F({{P}_{n}}\times {{K}_{m}})=m(8n-14){{(m-1)}^{3}}$,

(vi) $F({{C}_{n}}\times {{K}_{m}})=8nm{{(m-1)}^{3}}$.
\end{ex}

\subsection{Strong Product}
The strong product of two graphs $G_1$ and $G_2$ is denoted by ${{G}_{1}}\boxtimes {{G}_{2}}$. It has the vertex set $V(G_1)\times V(G_2)$ and any two vertices $({{u}_{p}},{{v}_{r}})$ and $({{u}_{q}},{{v}_{s}})$ are adjacent if and only if [${{u}_{p}}={{u}_{q}}$ and ${{v}_{r}}{{v}_{s}}\in E(G_2)$] or [${{v}_{r}}={{v}_{s}}$ and ${{u}_{p}}{{u}_{q}}\in E(G_1)$] or [${{u}_{p}}{{u}_{q}}\in E(G_1)$ and ${{v}_{r}}{{v}_{s}}\in E(G_2)$]. Note that if both $G_1$ and $G_2$ are connected then ${{G}_{1}}\boxtimes {{G}_{2}}$ is also connected. The degree of a vertex $(a,b)$ of ${{G}_{1}}\boxtimes {{G}_{2}}$ is given by
\[{{d}_{{{G}_{1}}\boxtimes {{G}_{2}}}}(a,b)={{d}_{{{G}_{1}}}}(a)+{{d}_{{{G}_{2}}}}(b)+{{d}_{{{G}_{1}}}}(a){{d}_{{{G}_{2}}}}(b).\]
In the following theorem we compute the F-index of the strong product of two graphs.
\begin{thm}

The F-index of ${{G}_{1}}\boxtimes {{G}_{2}}$ is given by
	\[F({{G}_{1}}\boxtimes {{G}_{2}})={{n}_{2}}F({{G}_{1}})+{{n}_{1}}F({{G}_{2}})+F({{G}_{1}})F({{G}_{2}})+6{{m}_{2}}{{M}_{1}}({{G}_{1}})+6{{m}_{1}}{{M}_{1}}({{G}_{2}})+6{{m}_{2}}F({{G}_{1}})\]
            \[+6{{m}_{1}}F({{G}_{2}})+3F({{G}_{2}}){{M}_{1}}({{G}_{1}})+3F({{G}_{1}}){{M}_{1}}({{G}_{2}})+6{{M}_{1}}({{G}_{1}}){{M}_{1}}({{G}_{2}}).\]
\end{thm}
\n\textit{Proof.} From definition of F-index, we have
\begin{eqnarray*}
F({{G}_{1}}\boxtimes {{G}_{2}})&=&\sum\limits_{({{v}_{1}},{{v}_{2}})\in V({{G}_{1}}\otimes {{G}_{2}})}{{{d}_{{{G}_{1}}\otimes {{G}_{2}}}}}{{({{v}_{1}},{{v}_{2}})}^{3}}\\
                     &=&\sum\limits_{v_1\in V(G_1)}\sum\limits_{v_2\in V(G_2)}[d_{G_1}(v_1)+d_{G_2}(v_2)+d_{G_1}(v_1)d_{G_2}(v_2)]^3\\
              &=&\sum\limits_{v_1\in V(G_1)}\sum\limits_{v_2\in V(G_2)} [d_{G_1}(v_1)^3+d_{G_2}              (v_2)^3+d_{G_1}(v_1)^3d_{G_2}(v_2)^3+3d_{G_1}(v_1)^2d_{G_2}(v_2)\\
                &&+3d_{G_1}(v_1)d_{G_2}(v_2)^2+3d_{G_1}(v_1)^3d_{G_2}(v_2)+3d_{G_1}(v_1)^3d_{G_2}(v_2)^2\\
            &&+3d_{G_2}(v_2)^3d_{G_1}(v_1)+3d_{G_1}(v_1)^2d_{G_2}(v_2)^3+6d_{G_1}(v_1)^2d_{G_2}(v_2)^2].
                \end{eqnarray*}
On simplification, we obtain the desired result.\qed

\subsection{Corona Product}
The corona product ${{G}_{1}}\odot {{G}_{2}}$ of two graphs $G_1$ and $G_2$ is obtained by taking one copy of ${{G}_{1}}$ and ${{n}_{1}}$ copies of ${{G}_{2}}$ and by joining each vertex of the $i$-th copy of ${{G}_{2}}$ to the $i$-th vertex of ${{G}_{1}}$, where $1\le i\le {{n}_{1}}$. The corona product of ${{G}_{1}}$ and ${{G}_{2}}$ has total $({{n}_{1}}{{n}_{2}}+{{n}_{1}})$ number of vertices and $({{m}_{1}}+{{n}_{1}}{{m}_{2}}+{{n}_{1}}{{n}_{2}})$ number of edges. Clearly, the corona product of two graphs is not commutative. Different topological indices under the corona product of graphs have already been studied by some researchers \cite{pat14,yar12}. It is easy to see that the degree of a vertex $v$ of ${{G}_{1}}\odot {{G}_{2}}$ is given by
\[{{d}_{{{G}_{1}}\odot{{G}_{2}}}}(v) = \left\{ \begin{array}{ll}
{{d}_{{{G}_{1}}}}(v)+{{n}_{2}},v\in V({{G}_{1}})\\[2mm]
{{d}_{{{G}_{2}}}}(v)+1,v\in V({{G}_{2,i}}), i=1,2,\ldots,{n_1}.
\end{array}\right.\]
where, ${{G}_{2,i}}$ is the i-th copy of the graph $G_2$.
In the following theorem, the F-index of the corona Product of two graphs is computed.
\begin{thm}

The F-index of ${{G}_{1}}\odot {{G}_{2}}$ is given by
\[F({{G}_{1}}\odot {{G}_{2}})=F({{G}_{1}})+{{n}_{1}}F({{G}_{2}})+3{{n}_{2}}{{M}_{1}}({{G}_{1}})+3{{n}_{1}}{{M}_{1}}({{G}_{2}})+6{{n}_{2}}^{2}{{m}_{1}}+6{{n}_{1}}{{m}_{2}}+{{n}_{1}}{{n}_{2}}({{n}_{2}}^{2}+1).\]
\end{thm}
\n\textit{Proof.} From definition of F-index, we have
\begin{eqnarray*}
F({{G}_{1}}\odot {{G}_{2}})&=&\sum\limits_{v\in V({{G}_{1}})}{{{({{d}_{{{G}_{1}}}}(v)+{{n}_{2}})}^{3}}}+{{n}_{1}}\sum\limits_{v\in V({{G}_{2}})}{{{({{d}_{{{G}_{2}}}}(v)+1)}^{3}}}\\
&=&\sum\limits_{v\in V({{G}_{1}})}{({{d}_{{{G}_{1}}}}{{(v)}^{3}}+3{{n}_{2}}{{d}_{{{G}_{1}}}}{{(v)}^{2}}+3{{n}_{2}}^{2}{{d}_{{{G}_{1}}}}(v)+{{n}_{2}}^{3})}\\
&&+{{n}_{1}}\sum\limits_{v\in V({{G}_{2}})}{({{d}_{{{G}_{2}}}}{{(v)}^{3}}+3{{d}_{{{G}_{2}}}}{{(v)}^{2}}+3{{d}_{{{G}_{2}}}}(v)+1)}\\
&=&\sum\limits_{v\in V({{G}_{1}})}{{{d}_{{{G}_{1}}}}{{(v)}^{3}}}+3{{n}_{2}}\sum\limits_{v\in V({{G}_{1}})}{{{d}_{{{G}_{1}}}}{{(v)}^{2}}}+3{{n}_{2}}^{2}\sum\limits_{v\in V({{G}_{1}})}{{{d}_{{{G}_{1}}}}(v)}+{{n}_{1}}{{{n}_{2}}^3}\\
&&+{{n}_{1}}\sum\limits_{v\in V({{G}_{1}})}{{{d}_{{{G}_{2}}}}{{(v)}^{3}}}+3{{n}_{1}}\sum\limits_{v\in V({{G}_{1}})}{{{d}_{{{G}_{2}}}}{{(v)}^{2}}}+3{{n}_{1}}\sum\limits_{v\in V({{G}_{1}})}{{{d}_{{{G}_{2}}}}(v)}+{{n}_{1}}{{n}_{2}}\\
&=&F({{G}_{1}})+3{{n}_{2}}{{M}_{1}}({{G}_{1}})+6{{n}_{2}}^{2}{{m}_{1}}+{{n}_{1}}{{n}_{2}}^{3}+{{n}_{1}}F({{G}_{2}})+3{{n}_{1}}{{M}_{1}}({{G}_{2}})\\
&&+6{{n}_{1}}{{m}_{2}}+{{n}_{1}}{{{n}_{2}}^3}+{{n}_{1}}{{n}_{2}}.
\end{eqnarray*}
   \qed

The $t$-thorny graph $G^t$ of a given graph $G$ is obtained by joining $t$-number of thorns (pendent edges) to each vertex of $G$. A variety of topological indices of thorn graphs have been studied by a number of researchers \cite{ali14,de12,nd12,de3}.It is well known that, the $t$-thorny graph of $G$ is defined as the corona product of $G$ and complement of complete graph with $t$ vertices $\bar{K_t}$. Thus from the previous theorem the following corollary follows.
\begin{cor}
The F-index of the $t$-thorny graph is given by
\[F({{G}^{t}})=F(G)+3t{{M}_{1}}(G)+6m{{t}^{2}}+n{{t}^{3}}+nt\]
where, $n$ and $m$ are number of vertices and edges of $G$, respectively.
\end{cor}

\begin{ex}
The F-index of $t$-thorny graph of ${{C}_{n}}$ and ${{P}_{n}}$ are given by

(i) $F({{C}_{n}}^{t})=n{{t}^{3}}+6n{{t}^{2}}+13nt+8n$,

(ii) $F({{P}_{n}}^{t})=n{{t}^{3}}+6n{{t}^{2}}-6{{t}^{2}}+13nt-18t+8n-14$.
\end{ex}

\begin{ex}
One of the hydrogen suppressed molecular graph is the bottleneck graph $(B)$  of a graph $G$, which is defined as the corona product of ${K}_{2}$ and $G$, where $G$ is a given graph. The F-index of bottleneck graph of $G$ is given by $F(B)=2F(G)+6{{M}_{1}}(G)+2n^3+6n^2+8n+12m+2$, where $n$ and $m$ are the number of vertices and edges of $G$, respectively.
\end{ex}

Next, as an application of corona product of graphs, we find the F index of some particular bridge graphs. Let $G_{1}, G_{2},...,G_{n}$ be a set of finite pairwise disjoint graphs. The bridge graph with respect to the vertices $v_{1}, v_{2},...,v_{n}$, denoted by $B(G_{1}, G_{2},...,G_{n};v_{1}, v_{2},...,v_{n})$ is the graph obtained by connecting the vertex $v_{i}$ of $G_{i}$ and the vertex $v_{i+1}$ of $G_{i+1}$ by an edge for all $i=1,2,...,n-1$. If $G_{1}\cong G_{2}\cong...\cong G_{n}$ and $v_{1}= v_{2}=...= v_{n}=v$, we define $G_{n}(G,v)=B(G, G,...,G;v, v,...,v)$. In particular, let $B_{n}=G_{n}(P_{3},v)$, where $v$ is the degree 2 vertex of $P_{3}$ and $T_{n,k}=G_{n}(C_{k},u)$  be special bridge graphs. Then from definition of corona product of graphs, $B_{n}\cong P_{n}\odot \overline{K}_{2}$ and $T_{n,3}\cong P_{n}\odot K_{2}$. Using the previous theorem the F index of these bridge graphs are calculated as follows.

\begin{ex}

(i) ${F(B_{m})}=66m-74$, for $m\geq2$

(ii) ${F(T_{m,3})}=80m-74$, for $m\geq2$

(iii) ${F(J_{n,m+1})}={n^3}m+6{n^2}m-6{n^2}+39nm+8m-18n-14$ , for $n\geq3$ and $m\geq2$.
\end{ex}

\subsection{Generalized Hierarchical Product}

As an extension of Cartesian product of graphs, Barriere et. al. introduced the generalized hierarchical product of graphs in 2009 \cite{barr09}. Several results on different topological indices under generalized hierarchical product of graphs are already studied \cite{elia13,arez13,luo14,nd14a}. Let $G_1$ and $G_2$ be two connected graphs and $\phi \ne U\subseteq V(G_2)$. Then the hierarchical product of $G_1$ and $G_2$ denoted by $G_1\Pi G_2(U)$, is the graph with vertex set $V(G_1)\times V(G_2)$ and any two vertices $(u,v)$ and $({u}',{v}')$ of $G_1(U)\Pi G_2$ are adjacent if and only if [$u={u}'\in V(G_1)$ and $v{v}'\in E(G_2)$] or [$v={v}'\in U$ and $u{u}'\in E(G_1)$].

Clearly, the degree of a vertex $(u_1,u_2)$ of ${{G}_{1}}\Pi {{G}_{2}(U)}$ is given by
\[{{d}_{G_1\Pi G_2(U)}}(u) = \left\{ \begin{array}{ll}
{{d}_{{{G}_{1}}}}(u_1)+{{d}_{{{G}_{2}}}}(u_2),{u_2}\in U\\[2mm]
{{d}_{{{G}_{2}}}}(u_2),{u_2}\in V({{G}_{2}})-U.
\end{array}\right.\]
In the following theorem we compute the F-index of the generalized hierarchical product of two graphs.
\begin{thm}
The F index of ${{G}_{1}}\Pi {{G}_{2}}(U)$ is given by

\[F({{G}_{1}}\Pi {{G}_{2}}(U))=|U|F({{G}_{1}})+{{n}_{1}}F({{G}_{2}})+3{{M}_{1}}({{G}_{1}})\sum\limits_{v\in U}{{{d}_{{{G}_{2}}}}(v)}+6{{m}_{1}}\sum\limits_{v\in U}{{{d}_{{{G}_{2}}}}{{(v)}^{2}}}.\]
\end{thm}

\n\textit{Proof.}
We have, from definition of F index
\begin{eqnarray*}
F({{G}_{1}}\Pi {{G}_{2}}(U))&=&\sum\limits_{({{v}_{1}},{{v}_{2}})\in V({{G}_{1}}\Pi {{G}_{2}}(U))}{{{d}_{{{G}_{1}}\Pi {{G}_{2}}(U)}}}{{({{v}_{1}},{{v}_{2}})}^{3}}\\
&=&\sum\limits_{{{v}_{1}}\in V({{G}_{1}})}{\sum\limits_{{{v}_{2}}\in U}{{{[{{d}_{{{G}_{1}}}}({{v}_{1}})+{{d}_{{{G}_{2}}}}({{v}_{2}})]}^{3}}}}+\sum\limits_{{{v}_{1}}\in V({{G}_{1}})}{\sum\limits_{{{v}_{2}}\in V({{G}_{2}})-U}{{{d}_{{{G}_{2}}}}{{({{v}_{2}})}^{3}}}}\\
&=&\sum\limits_{{{v}_{1}}\in V({{G}_{1}})}{\sum\limits_{{{v}_{2}}\in U}{[{{d}_{{{G}_{1}}}}{{({{v}_{1}})}^{3}}+{{d}_{{{G}_{2}}}}{{({{v}_{2}})}^{3}}+3{{d}_{{{G}_{1}}}}{{({{v}_{1}})}^{2}}{{d}_{{{G}_{2}}}}({{v}_{2}})+3{{d}_{{{G}_{1}}}}({{v}_{1}}){{d}_{{{G}_{2}}}}{{({{v}_{2}})}^{2}}]}}\\
&&+\sum\limits_{{{v}_{1}}\in V({{G}_{1}})}{\sum\limits_{{{v}_{2}}\in V({{G}_{2}})-U}{{{d}_{{{G}_{2}}}}{{({{v}_{2}})}^{3}}}}\\
&=&|U|F({{G}_{1}})+{{n}_{1}}F({{G}_{2}})+3{{M}_{1}}({{G}_{1}})\sum\limits_{v\in U}{{{d}_{{{G}_{2}}}}(v)}+6{{m}_{1}}\sum\limits_{v\in U}{{{d}_{{{G}_{2}}}}{{(v)}^{2}}},\end{eqnarray*}
which completes the proof.   \qed

From definition of the Cartesian product of graphs, it is clear that the Cartesian product of graphs is a special case of generalized hierarchical product of graphs, that is, if $U=V(G_2)$, then $G_1\Pi G_2(U)\cong G_1\otimes G_2$. So from the previous theorem we can also obtain the Theorem 3, for $k=2$.

The cluster product of two graphs $G_1$ and $G_2$, denoted by $G_1\left\{ G_2 \right\}$, is obtained by taking one copy of $G_1$ and $n_1$ copies of a rooted graph $G_2$ and by identifying the root of the $i$-th copy of $G_2$ with the $i$-th vertex of $G_1$, $i=1,2,...,n_1$. From definition of cluster product of graphs, $|V(G_1\left\{ G_2 \right\})|=n_1{n_2}$ and $|E(G_1\left\{ G_2 \right\})|=(m_2+n_1{m_2})$. Suppose the root vertex of $G_2$ is denoted by $x$. Note that, if $U=\left\{ x \right\}\subset V(G_1)$ then $G_1\left\{ G_2 \right\}\cong G_1\Pi G_2(U)\cong G_1\Pi G_2(\left\{ x \right\})$. Then from the previous theorem the following results follow.

\begin{cor}
Let ${{G}_{1}}$ and ${{G}_{2}}$ be two connected graphs and $x$ be the root vertex of ${{G}_{2}}$, then
\[F({{G}_{1}}\{{{G}_{2}}\})=F({{G}_{1}})+{{n}_{1}}F({{G}_{2}})+3{{M}_{1}}({{G}_{1}}){{d}_{{{G}_{2}}}}(x)+6{{m}_{1}}{{d}_{{{G}_{2}}}}{{(x)}^{2}}.\]
\end{cor}

Note that, if $U=\left\{ x \right\}$, $x$ is the root vertex of the graph $G_2$, then ${{G}_{1}}\Pi {{G}_{2}}(U)={{G}_{1}}\Pi {{G}_{2}}$, the (standard) hierarchical product of two graphs. It is easy to see that ${{G}_{1}}\Pi {{G}_{2}}\cong{{G}_{1}}\{{{G}_{2}}\}.$
\begin{ex}
The square comb lattice $Cq(N)$ with $N={{n}^{2}}$ vertices can be represented as the cluster product ${{P}_{n}}\{{{P}_{n}}\}$, where the root of ${{P}_{n}}$ is one of its vertices of degree 1. Then from the previous corollary, the F-index of $Cq(N)$ is given by
$F(Cq(N))=8{{n}^{2}}+12n-38.$
\end{ex}
\begin{ex}
The sun graph is defined as ${Sun_{(m,n)}}={{C}_{m}}\{{{P}_{n+1}}\}$, such that ${P_{n+1}}$ is rooted at a vertex of degree one. Then the F-index of ${Sun_{(m,n)}}$ is given by ${Sun_{(m,n)}}=4m(2n+5)$.
\end{ex}

\subsection{Disjunction}
The disjunction of two graphs ${{G}_{1}}$ and ${{G}_{2}}$, denoted by ${{G}_{1}}\wedge{{G}_{2}}$, consists of the vertex set $V(G_1)\times V(G_2)$ and two vertices $(u_1,v_1)$ and $(u_2,v_2)$ are adjacent whenever ${u_1}{u_2}\in{E(G_1)}$ or ${v_1}{v_2}\in{E(G_2)}$. Clearly, the degree of a vertex $(u_1,u_2)$ of ${{G}_{1}}\wedge {{G}_{2}}$ is given by
\[{{d}_{{{G}_{1}}\wedge {{G}_{2}}}}(u_1,u_2)={n_1}{{d}_{{{G}_{1}}}}(u_1)+{n_2}{{d}_{{{G}_{2}}}}(u_2)-{{d}_{{{G}_{1}}}}(u_1){{d}_{{{G}_{2}}}}(u_2).\]
In the following theorem we obtain the F-index of the disjunction of two graphs.

\begin{thm} The F-index of ${{G}_{1}}\wedge {{G}_{2}}$ is given by
\begin{eqnarray*}
F({{G}_{1}}\wedge {{G}_{2}})&=&{{n}_{2}}^{4}F({{G}_{1}})+{{n}_{1}}^{4}F({{G}_{2}})-F({{G}_{1}})F({{G}_{2}})+6{{n}_{1}}{{n}_{2}}^{2}{{m}_{2}}{{M}_{1}}({{G}_{1}})+6{{n}_{1}}^{2}{{n}_{2}}{{m}_{1}}{{M}_{1}}({{G}_{2}})\\
&&+3{{n}_{2}}F({{G}_{1}}){{M}_{1}}({{G}_{2}})+3{{n}_{1}}F({{G}_{2}}){{M}_{1}}({{G}_{1}})-6{{n}_{2}}^{2}{{m}_{2}}F({{G}_{1}})-6{{n}_{1}}^{2}{{m}_{1}}F({{G}_{2}})\\
&&-6{{n}_{1}}{{n}_{2}}{{M}_{1}}({{G}_{1}}){{M}_{1}}({{G}_{2}}).\\
\end{eqnarray*}
\end{thm}

\n\textit{Proof.}
We have, from definition of F-index
\begin{eqnarray*}
F({{G}_{1}}\wedge {{G}_{2}})&=&\sum\limits_{({{u}_{1}},{{u}_{2}})\in V({{G}_{1}}\wedge {{G}_{2}})}{{{d}_{{{G}_{1}}\wedge{{G}_{2}}}}{{({{u}_{1}},{{u}_{2}})}^{3}}}\\
&=&\sum\limits_{{{u}_{1}}\in V({{G}_{1}})}{\sum\limits_{{{u}_{2}}\in V({{G}_{2}})}{({{n}_{2}}{{d}_{{{G}_{1}}}}({{u}_{1}})+}}{{n}_{1}}{{d}_{{{G}_{2}}}}({{u}_{2}})-{{d}_{{{G}_{1}}}}({{u}_{1}}){{d}_{{{G}_{2}}}}({{u}_{2}}){{)}^{3}}\\
&=&\sum\limits_{{{u}_{1}}\in V({{G}_{1}})}{\sum\limits_{{{u}_{2}}\in V({{G}_{2}})}{({{n}_{2}}^{3}{{d}_{{{G}_{1}}}}{{({{u}_{1}})}^{3}}+}}{{n}_{1}}^{3}{{d}_{{{G}_{2}}}}{{({{u}_{2}})}^{3}}-{{d}_{{{G}_{1}}}}{{({{u}_{1}})}^{3}}{{d}_{{{G}_{2}}}}{{({{u}_{2}})}^{3}}\\
&&+3{{n}_{1}}{{n}_{2}}^{2}{{d}_{{{G}_{1}}}}{{({{u}_{1}})}^{2}}{{d}_{{{G}_{2}}}}({{u}_{2}})+3{{n}_{1}}^{2}{{n}_{2}}{{d}_{{{G}_{1}}}}({{u}_{1}}){{d}_{{{G}_{2}}}}{{({{u}_{2}})}^{2}}-3{{n}_{2}}^{2}{{d}_{{{G}_{1}}}}{{({{u}_{1}})}^{3}}{{d}_{{{G}_{2}}}}({{u}_{2}})\\
&&-3{{n}_{1}}^{2}{{d}_{{{G}_{1}}}}({{u}_{1}}){{d}_{{{G}_{2}}}}{{({{u}_{2}})}^{3}}+3{{n}_{2}}{{d}_{{{G}_{1}}}}{{({{u}_{1}})}^{3}}{{d}_{{{G}_{2}}}}{{({{u}_{2}})}^{2}}+3{{n}_{1}}{{d}_{{{G}_{1}}}}{{({{u}_{1}})}^{2}}{{d}_{{{G}_{2}}}}{{({{u}_{2}})}^{3}}\\
&&-6{{n}_{1}}{{n}_{2}}{{d}_{{{G}_{1}}}}{{({{u}_{1}})}^{2}}{{d}_{{{G}_{2}}}}{{({{u}_{2}})}^{2}}).
\end{eqnarray*}

After simple calculations, we obtain the desired result.\qed

\subsection{Symmetric Difference}

The symmetric difference, of two graphs ${{G}_{1}}$ and ${{G}_{2}}$ is denoted by
${{G}_{1}}\oplus {{G}_{2}}$, so that $|V({{G}_{1}}\oplus {{G}_{2}})|=|V({{G}_{1}})|\times |V({{G}_{2}})|$ and
	\[E({{G}_{1}}\oplus {{G}_{2}})=\{({{u}_{1}},{{u}_{2}})({{v}_{1}},{{v}_{2}}):{{u}_{1}}{{v}_{1}}\in E({{G}_{1}})\,\mbox{ or } {{u}_{2}}{{v}_{2}}\in E({{G}_{2}})\mbox{ but not both}\}.\]
From definition of symmetric difference it is clear that
\[{{d}_{{{G}_{1}}\oplus {{G}_{2}}}}({{v}_{1}},{{v}_{2}})={{n}_{2}}{{d}_{{{G}_{1}}}}({{v}_{1}})+{{n}_{1}}{{d}_{{{G}_{2}}}}({{v}_{2}})-2{{d}_{{{G}_{1}}}}({{v}_{1}}){{d}_{{{G}_{2}}}}({{v}_{2}}).\]
In the following theorem we obtain the F-index of the symmetric difference of two graphs.
\begin{thm} The F-index of ${{G}_{1}}\oplus {{G}_{2}}$ is given by
\begin{eqnarray*}
F({{G}_{1}}\oplus {{G}_{2}})&=&{{n}_{2}}^{4}F({{G}_{1}})+{{n}_{1}}^{4}F({{G}_{2}})-8F({{G}_{1}})F({{G}_{2}})+6{{n}_{1}}{{n}_{2}}^{2}{{m}_{2}}{{M}_{1}}({{G}_{1}})+6{{n}_{1}}^{2}{{n}_{2}}{{m}_{1}}{{M}_{1}}({{G}_{2}})\\
&&+12{{n}_{2}}F({{G}_{1}}){{M}_{1}}({{G}_{2}})+12{{n}_{1}}F({{G}_{2}}){{M}_{1}}({{G}_{1}})-12{{n}_{2}}^{2}{{m}_{2}}F({{G}_{1}})-12{{n}_{1}}^{2}{{m}_{1}}F({{G}_{2}})\\
&&-12{{n}_{1}}{{n}_{2}}{{M}_{1}}({{G}_{1}}){{M}_{1}}({{G}_{2}}).
\end{eqnarray*}
\end{thm}
\n\textit{Proof.} By definition of F-index, we have
\begin{eqnarray*}
F({{G}_{1}}\oplus {{G}_{2}})&=&\sum\limits_{({{u}_{1}},{{u}_{2}})\in V({{G}_{1}}\wedge {{G}_{2}})}{{{d}_{{{G}_{1}}\oplus {{G}_{2}}}}{{({{u}_{1}},{{u}_{2}})}^{3}}}\\
&=&\sum\limits_{{{u}_{1}}\in V({{G}_{1}})}{\sum\limits_{{{u}_{2}}\in V({{G}_{2}})}{({{n}_{2}}{{d}_{{{G}_{1}}}}({{u}_{1}})+}}{{n}_{1}}{{d}_{{{G}_{2}}}}({{u}_{2}})-2{{d}_{{{G}_{1}}}}({{u}_{1}}){{d}_{{{G}_{2}}}}({{u}_{2}}){{)}^{3}}\\
&=&\sum\limits_{{{u}_{1}}\in V({{G}_{1}})}{\sum\limits_{{{u}_{2}}\in V({{G}_{2}})}{({{n}_{2}}^{3}{{d}_{{{G}_{1}}}}{{({{u}_{1}})}^{3}}+}}{{n}_{1}}^{3}{{d}_{{{G}_{2}}}}{{({{u}_{2}})}^{3}}-8{{d}_{{{G}_{1}}}}{{({{u}_{1}})}^{3}}{{d}_{{{G}_{2}}}}{{({{u}_{2}})}^{3}}\\
&&+3{{n}_{1}}{{n}_{2}}^{2}{{d}_{{{G}_{1}}}}{{({{u}_{1}})}^{2}}{{d}_{{{G}_{2}}}}({{u}_{2}})+3{{n}_{1}}^{2}{{n}_{2}}{{d}_{{{G}_{1}}}}({{u}_{1}}){{d}_{{{G}_{2}}}}{{({{u}_{2}})}^{2}}-6{{n}_{2}}^{2}{{d}_{{{G}_{1}}}}{{({{u}_{1}})}^{3}}{{d}_{{{G}_{2}}}}({{u}_{2}})\\ &&-6{{n}_{1}}^{2}{{d}_{{{G}_{1}}}}({{u}_{1}}){{d}_{{{G}_{2}}}}{{({{u}_{2}})}^{3}}+12{{n}_{2}}{{d}_{{{G}_{1}}}}{{({{u}_{1}})}^{3}}{{d}_{{{G}_{2}}}}{{({{u}_{2}})}^{2}}+12{{n}_{1}}{{d}_{{{G}_{1}}}}{{({{u}_{1}})}^{2}}{{d}_{{{G}_{2}}}}{{({{u}_{2}})}^{3}}\\
&&-12{{n}_{1}}{{n}_{2}}{{d}_{{{G}_{1}}}}{{({{u}_{1}})}^{2}}{{d}_{{{G}_{2}}}}{{({{u}_{2}})}^{2}}),
\end{eqnarray*}
from where we obtain the desired result.\qed

\subsection{Splice and Link of Graphs}

Let ${{v}_{1}}\in {{V}({G}_{1})}$ and ${{v}_{2}}\in {{V}({G}_{2})}$ be two given vertices of ${{G}_{1}}$ and ${{G}_{2}}$, respectively. A splice of ${{G}_{1}}$ and ${{G}_{2}}$ at the vertices ${{v}_{1}}$ and ${{v}_{2}}$, denoted by $({{G}_{1}}\bullet {{G}_{2}})({{v}_{1}},{{v}_{2}})$, was introduced by Dosli\'{c} \cite{dosl05}, and is obtained by identifying the vertices ${{v}_{1}}$ and ${{v}_{2}}$ in the union of ${{G}_{1}}$ and ${{G}_{2}}$. The vertex set of $({{G}_{1}}\bullet {{G}_{2}})({{v}_{1}},{{v}_{2}})$ is given by
$V(({{G}_{1}}\bullet {{G}_{2}})({{v}_{1}},{{v}_{2}}))=[V({{G}_{1}})\backslash \{{{v}_{1}}\}]\cup [V({{G}_{2}})\backslash \{{{v}_{2}}\}]\cup \{{{v}_{12}}\}$, where we denote the vertex obtained by identifying ${{v}_{1}}$ and ${{v}_{2}}$ by ${{v}_{12}}$.  From the construction of splice graphs it is clear that
\[
{d_{({{G_1} \bullet {G_2}})({v_1},{v_2})}(v)} = \left\{ \begin{array}{ll}
{{d_{{G_i}}}(v)},&\mbox{for } v \in V({G_i})\mbox{ and } v \ne {v_i},\\
{{d_{{G_1}}}({v_1}) + {d_{{G_2}}}({v_2})},&\mbox{for } v = {v_{12}}.
\end{array} \right.
\]

Similarly a link of ${{G}_{1}}$ and ${{G}_{2}}$ at the vertices ${{v}_{1}}$ and ${{v}_{2}}$ is denoted by $({{G}_{1}}\sim {{G}_{2}})({{v}_{1}},{{v}_{2}})$ and is obtained by joining the vertices ${{v}_{1}}$ and ${{v}_{2}}$ in the union of ${{G}_{1}}$ and ${{G}_{2}}$. From the construction of link graphs it is clear that
\[d_{(G_1\sim G_2){(v_1,v_2)}}(v)=\left\{\begin{array}{ll}d_{G_i}(v),& v\in V(G_i)\mbox{ and } v\neq v_i, i=1,2,\\ d_{G_i}(v)+1,& v = v_i, i=1,2.\end{array}\right.\]
In the following, we find the F-index of splice and link of two graphs ${{G}_{1}}$ and ${{G}_{2}}$ at the vertices $v_1$ and $z$.

\begin{thm}
 The F-index of $({{G}_{1}}\bullet {{G}_{2}})(v_1,v_2)$ is given by
\[F(({{G}_{1}}\bullet {\ }{{G}_{2}})(v_1,v_2))=F({{G}_{1}})+F({{G}_{2}})+3{{d}_{{{G}_{1}}}}(v_1){{d}_{{{G}_{2}}}}(v_2)({{d}_{{{G}_{1}}}}(v_1)+{{d}_{{{G}_{2}}}}(v_2)).\]
\end{thm}
\n\textit{Proof.} From the definition of F-index we have
\begin{eqnarray*}
F(({{G}_{1}}\bullet {{G}_{2}})(v_1,v_2))&=&{{\sum\limits_{v\in V({{G}_{1}}),v\ne v_1}{{{d}_{{{G}_{1}}}}(v)}}^{3}}+{{\sum\limits_{v\in V({{G}_{2}}),v\ne v_2}{{{d}_{{{G}_{2}}}}(v)}}^{3}}+{{({{d}_{{{G}_{1}}}}(v_1)+{{d}_{{{G}_{2}}}}(v_2))}^{3}}\\
&=&{{\sum\limits_{v\in V({{G}_{1}})}{{{d}_{{{G}_{1}}}}(v)}}^{3}}+{{\sum\limits_{v\in V({{G}_{2}})}{{{d}_{{{G}_{2}}}}(v)}}^{3}}+3{{d}_{{{G}_{1}}}}{{(v_1)}^{2}}{{d}_{{{G}_{2}}}}(v_2)+3{{d}_{{{G}_{1}}}}(v_1){{d}_{{{G}_{2}}}}{{(v_2)}^{2}}.
\end{eqnarray*}
From above we obtain the desired result after simple calculations.\qed

\begin{thm} The F-index of $({{G}_{1}}\sim{\ }{{G}_{2}})(v_1,v_2)$ is given by
	
\[F(({{G}_{1}}\sim{\ }{{G}_{2}})(v_1,v_2))= F({{G}_{1}})+F({{G}_{2}})+3({{d}_{{{G}_{1}}}}(v_1)+{{d}_{{{G}_{2}}}}(v_2))+3({{d}_{{{G}_{1}}}}{{(v_1)}^{2}}+{{d}_{{{G}_{2}}}}{{(v_2)}^{2}})+2.\]
\end{thm}
\n\textit{Proof.} From the definition of F-index, we have
\begin{eqnarray*}
F(({{G}_{1}}\sim {{G}_{2}})(v_1,v_2))&=&{{\sum\limits_{v\in V({{G}_{1}}),v\ne v_1}{{{d}_{{{G}_{1}}}}(v)}}^{3}}+{{\sum\limits_{v\in V({{G}_{2}}),v\ne v_2}{{{d}_{{{G}_{2}}}}(v)}}^{3}}+{{({{d}_{{{G}_{1}}}}(v_1)+1)}^{3}}+{{({{d}_{{{G}_{2}}}}(v_2)+1)}^{3}}\\
&=&F({{G}_{1}})+F({{G}_{2}})+3{{d}_{{{G}_{1}}}}(v_1)+3{{d}_{{{G}_{1}}}}{{(v_1)}^{2}}+1+3{{d}_{{{G}_{2}}}}(v_2)+3{{d}_{{{G}_{2}}}}{{(v_2)}^{2}}+1,
\end{eqnarray*}
from where the desired result follows. \qed

\section{Conclusion}
In this paper, we derive some explicit expression of the F-index of different graph operations such as union, join, Cartesian product, composition, tensor product, strong product, corona product, generalized hierarchical product, disjunction, symmetric difference, splice and link of graphs. Also we apply our results to compute the F-index for some important classes of molecular graphs and nano-structures. For further study, F-index of some other graph operations and for different composite graphs can be computed.

\end{document}